\documentclass[11pt]{article}
\usepackage[cp850]{inputenc}
\usepackage{amssymb}
\usepackage{latexsym}
\def\Z{\mathbb{Z}}
\def\N{\mathbb{N}}
\def\R{\mathbb{R}}
\def\C{\mathbb{C}}
\def\s{\frak{s}}

\begin{document}
\setlength{\parindent}{0pt}
\setlength{\parskip}{0.4cm}

\begin{center}

\large{\bf Counting Lattice Points by means of the Residue Theorem}
\footnote{This paper appeared in {\it Ramanujan Journal} {\bf 4}, no.~3 (2000), 299--310. 
          Parts of this work appeared in the author's Ph.D. thesis, written under supervision of Sinai Robins.} 

\normalsize{\sc Matthias Beck}

\end{center}

\footnotesize {\bf Abstract.}  
We use the residue theorem to derive an expression for the number of lattice points 
in a dilated $n$-dimensional tetrahedron with vertices 
at lattice points on each coordinate axis and the origin. 
This expression is known as the Ehrhart polynomial. 
We show that it is a polynomial in $t$, where $t$ is the integral dilation parameter. 
We prove the 
Ehrhart-Macdonald reciprocity law for these tetrahedra, relating the Ehrhart 
polynomials of the interior and the closure of the tetrahedra. To illustrate 
our method, we compute the Ehrhart coefficient for codimension 2. Finally, 
we show how our ideas can be used to compute the Ehrhart polynomial for an 
arbitrary convex lattice polytope.

\normalsize

\vspace{1cm}

{\bf 1. Introduction}

Let $\Z^{ n } \subset \R^{ n } $ be the $n$-dimensional integer lattice, and 
\[ {\cal P} = \left\{ (x_{ 1 } , \dots , x_{ n }) \in \R^{ n }  : \sum_{ k=1 }^{ n } \frac{ x_{ k } }{ a_{ k } } < 1 \ \mbox{and all} \ x_{ k } > 0 \right\} \]
the $n$-dimensional open tetrahedron with vertices 
$ (0, \dots , 0) $, $ \left( a_{ 1 } , 0 , \dots , 0 \right)$, $ \left( 0, a_{ 2 } , 0, \dots , 0 \right) , \dots , \left( 0, \dots , 0, a_{ n }  \right) $. 
where $a_{ 1 } ,$ $ \dots ,$ $a_{ n } $ are positive integers. 
For $t \in \N$, denote by $ L ( {\cal P} , t ) $ and $ L ( \overline{\cal P} , t ) $ the number of lattice 
points in the dilated polytope $ t {\cal P} $ and its closure, respectively. 

Ehrhart (\cite{ehrhart}) proved (actually for a general lattice polytope) that both $ L ( {\cal P} , t ) $ 
and $ L ( \overline{\cal P} , t ) $ are polynomials in $t$ of degree $n$. 
Moreover, he determined the two leading coefficients and the constant. 
The leading coefficient is Vol(${\cal P}$), and the second coefficient is 1/2 
Vol($\partial {\cal P}$), which is half the surface area of ${\cal P}$ normalized with respect to 
the sublattice on each face of ${\cal P}$. The constant coefficient equals $\chi ( {\cal P}) $, 
the Euler characteristic of ${\cal P}$. 
The other coefficients of $ L ( {\cal P} , t ) $ and $ L ( \overline{\cal P} , t ) $ are not as 
easily accessible. In fact, a method of computing these coefficients was unknown until 
quite recently (\cite{barv}, \cite{brion}, \cite{robins}, \cite{kantor}, \cite{puk}). 
In this paper, we present an elementary method for computing the Ehrhart polynomials of 
${\cal P}$ and $\overline{\cal P}$ using the residue theorem.
We verify the Ehrhart-Macdonald reciprocity law for these $n$-dimensional tetrahedra. To
illustrate our method, we compute the first nontrivial coefficient, $c_{ n-2 } $, of the Ehrhart 
polynomial. Finally, we show how our ideas can be used to compute the Ehrhart polynomial for an 
arbitrary convex lattice polytope.

\vspace{1cm}

{\bf 2. The main idea}

Let's start with $L ( \overline{\cal P} , t )$; that is, we consider the closure of our dilated tetrahedron $t {\cal P}$. 
We introduce the notation
\[ A := a_{1} \cdots a_{n} , \quad A_{k} := a_{1} \cdots {\hat a}_{k} \cdots a_{n} , \]
where ${\hat a}_{k}$ means we omit the factor $a_{k} $. Then we can write 
\begin{eqnarray*} &\mbox{}& L ( \overline{\cal P} , t ) = \mbox{card} \left\{ (m_{ 1 } , \dots , m_{ n }) \in \Z^{ n }  : \sum_{ k=1 }^{ n } \frac{ m_{ k } }{ a_{ k } } \leq t \ \mbox{and all} \ m_{ k } \geq 0 \right\} \\
                  &\mbox{}& \quad = \mbox{card}  \left\{ (m_{ 1 } , \dots , m_{ n }, m) \in \Z^{ n+1 }  : \begin{array}{l} \sum_{ k=1 }^{ n } m_{ k } A_{k} + m = t A \\ m_{ k }, m \geq 0 \end{array} \right\} . \end{eqnarray*} 
We can interpret $L ( \overline{\cal P} , t )$ as the Taylor coefficient of $z^{ t A }$ for the function
\begin{eqnarray*} &\mbox{}& \left( 1 + z^{ A_{ 1 } } + z^{ 2 A_{ 1 } } + \dots  \right) \left( 1 + z^{ A_{ 2 } } + z^{ 2 A_{ 2 } } + \dots  \right) \cdots \\ 
                  &\mbox{}& \quad \cdot \left( 1 + z^{ A_{ n } } + z^{ 2 A_{ n } } + \dots  \right) \left( 1 + z + z^{ 2 } + \dots \right) \\
                  &\mbox{}& = \frac{ 1 }{ 1 - z^{ A_{ 1 } } } \frac{ 1 }{ 1 - z^{ A_{ 2 } } } \cdots \frac{ 1 }{ 1 - z^{ A_{ n } } } \frac{ 1 }{ 1 - z } \ . \end{eqnarray*} 
Equivalently, 
\[ L ( \overline{\cal P} , t ) = \mbox{Res} \left( \frac{ z^{ - t A - 1 } }{ \left( 1 - z^{ A_{ 1 } }  \right) \left( 1 - z^{ A_{ 2 } }  \right) \cdots \left( 1 - z^{ A_{ n } }  \right) \left( 1 - z \right) } , z=0 \right) . \]
To reduce the number of poles, it is convenient to change this function slightly; this residue is clearly equal to 
\[ \mbox{Res} \left( \frac{ z^{ - t A } - 1 }{ \left( 1 - z^{ A_{ 1 } }  \right) \left( 1 - z^{ A_{ 2 } }  \right) \cdots \left( 1 - z^{ A_{ n } }  \right) \left( 1 - z \right) z }  , z=0 \right) + 1 . \]
If this expression counts the number of lattice points in $\overline{t \cal P}$, then all we have to do is compute the other residues of 
\[ f_{ -t } (z) := \frac{ z^{ - t A } - 1 }{ \left( 1 - z^{ A_{ 1 } }  \right) \left( 1 - z^{ A_{ 2 } }  \right) \cdots \left( 1 - z^{ A_{ n } }  \right) \left( 1 - z \right) z } \]
and use the residue theorem for the sphere $\C \cup \{ \infty \}$. 
In this notation, 
\begin{equation}\label{closure} L ( \overline{\cal P} , t ) = \mbox{Res} \left( f_{ -t } (z) , z=0 \right) + 1 \ . \end{equation}
The only poles of $f_{ -t } $ are at 0, 1 and the roots of unity in 
\[ \Omega := \left\{ z \in \C \setminus \{ 1 \} : z^{\frac{ A }{ a_{ k } a_{ j }  } } = 1 , 1 \leq k < j \leq n \right\} \ . \] 
Note that Res($f_{-t}, z=\infty$) = 0, so that the residue theorem gives us the first half of our main result:

{\bf Theorem 1(a).} 
\[ L ( \overline{\cal P} , t ) = 1 - \mbox{Res} \left( f_{ -t } (z) , z=1 \right) - \sum_{ \lambda \in \Omega } \mbox{Res} \left( f_{ -t } (z) , z=\lambda \right) \] 

{\it Remarks.} 1. The residue at z=1 can be calculated easily: 
\begin{eqnarray*} &\mbox{}& \mbox{Res} \left( f_{ -t } (z) , z=1 \right) = \mbox{Res} \left( e^{ z }  f_{ -t } (e^{ z } ) , z=0 \right)  \\
                  &\mbox{}& \quad = \mbox{Res} \left( \frac{ e^{ -t A z } - 1 }{ \left( 1 - e^{ A_{ 1 } z } \right) \left( 1 - e^{ A_{ 2 } z } \right) \cdots \left( 1 - e^{ A_{ n } z } \right) \left( 1 - e^{ z } \right) } , z=0 \right) . \end{eqnarray*} 
To facilitate the computation in higher dimensions, one can use mathematics software such as {\tt Maple} or {\tt Mathematica}.
It is easy to see that Res$ \left( f_{ -t } (z) , z=1 \right) $ is a polynomial in $t$ whose coefficients 
are rational expressions in $ a_{ 1 } , \dots, a_{ n } $. 

2. The residues at the roots of unity in $\Omega$ are in general not as easy to compute. 
They give rise to Dedekind-like sums and their higher dimensional analogues, as we will 
illustrate in section 4. There is, however, one feature we can read off from these residues immediately, 
the dependency on the dilation parameter $t$: 

{\bf Corollary 2 (Ehrhart)}. {\it $ L ( \overline{\cal P} , t ) $ is a polynomial in $t$.}

With Corollary 3 below, this will also imply that $ L ( {\cal P} , t ) $ is a polynomial.

{\it Proof.} 
Let $ \lambda \in \Omega $ be a $B$'th root of unity, where $B$ is the product of some of the $a_{ k }$. 
Now express $ z^{-tA} $ in terms 
of its power series about $ z=\lambda $. The coefficients of this power series involve 
various derivatives of $ z^{-tA} $, evaluated at $ z=\lambda $. Here we can introduce a change 
of variable: $ z = w^{ \frac{ 1 }{ B }  } = \exp \left( \frac{ 1 }{ B } \log w \right) $, 
where we choose a suitable branch of the logarithm such that 
$ \exp \left( \frac{ 1 }{ B } \log (1) \right) = \lambda$. The terms 
depending on $t$ in the power series of $ z^{-tA} $ consist therefore of derivatives of 
the function $ z^{-tA/B} $, evaluated at $z=1$. From this it is easy to 
see that the coefficients of the power series of $ z^{-tA} $ are polynomials in $t$. 
The fact that $ L ( \overline{\cal P} , t ) $ is simply the sum of all these residues, 
finally, gives the statement. \hfill {} $\Box$

For the computation of $L ( {\cal P} , t )$ (the number of lattice points in the {\it interior} 
of our tetrahedron ${t \cal P}$), we similarly write 
\begin{eqnarray*} &\mbox{}& L ( {\cal P} , t ) = \mbox{card} \left\{ (m_{ 1 } , \dots , m_{ n }) \in \Z^{ n }  : \sum_{ k=1 }^{ n } \frac{ m_{ k } }{ a_{ k } } < t \ \mbox{and all} \ m_{ k } > 0 \right\} \\
                  &\mbox{}& \qquad = \mbox{card} \left\{ (m_{ 1 } , \dots , m_{ n }, m) \in \Z^{ n+1 }  : \begin{array}{l} \sum_{ k=1 }^{ n } m_{ k } A_{k} + m = t A \\ m_{ k }, m > 0 \end{array} \right\} . \end{eqnarray*} 
Now we can interpret $L ( {\cal P} , t )$ as the Taylor coefficient of $z^{ t A }$ for the function
\begin{eqnarray*} &\mbox{}& \left( z^{ A_{ 1 } } + z^{ 2 A_{ 1 } } + \dots  \right) \left( z^{ A_{ 2 } } + z^{ 2 A_{ 2 } } + \dots  \right) \cdots \left( z^{ A_{ n } } + z^{ 2 A_{ n } } + \dots  \right) \left( z + z^{ 2 } + \dots \right) \\
                  &\mbox{}& \quad = \frac{ z^{ A_{ 1 } } }{ 1 - z^{ A_{ 1 } } } \ \frac{ z^{ A_{ 2 } } }{ 1 - z^{ A_{ 2 } } } \cdots \frac{ z^{ A_{ n } } }{ 1 - z^{ A_{ n } } } \ \frac{ z }{ 1 - z } \ , \end{eqnarray*} 
or equivalently as 
\begin{eqnarray*} &\mbox{}& \mbox{Res} \left( \frac{ z^{ A_{ 1 } } }{ 1 - z^{ A_{ 1 } } } \ \frac{ z^{ A_{ 2 } } }{ 1 - z^{ A_{ 2 } } } \cdots \frac{ z^{ A_{ n } } }{ 1 - z^{ A_{ n } } } \ \frac{ z }{ 1 - z } \ z^{ - t A - 1 } , z=0 \right)  \\
                  &\mbox{}& = \mbox{Res} \left( \frac{ z^{ A_{ 1 } } }{ 1 - z^{ A_{ 1 } } } \ \frac{ z^{ A_{ 2 } } }{ 1 - z^{ A_{ 2 } } } \cdots \frac{ z^{ A_{ n } } }{ 1 - z^{ A_{ n } } } \ \frac{ z }{ 1 - z } \ \frac{ z^{ - t A } - 1 }{ z } , z=0 \right) \\
                  &\mbox{}& = \mbox{Res} \left( \frac{ -1 }{ z^{ 2 }  } \ \frac{ 1 }{ z^{ A_{ 1 } } - 1 } \ \frac{ 1 }{ z^{ A_{ 2 } } - 1 } \cdots \frac{ 1 }{ z^{ A_{ n } } - 1 } \ \frac{ 1 }{ z - 1 } \ z \left( z^{ t A } - 1 \right) , z=\infty \right) . \end{eqnarray*} 
To be able to use the residue theorem, this time we have to consider the function
\[ - \frac{ 1 }{ z^{ A_{ 1 } } - 1 } \ \frac{ 1 }{ z^{ A_{ 2 } } - 1 } \cdots \frac{ 1 }{ z^{ A_{ n } } - 1 } \ \frac{ 1 }{ z - 1 } \ \frac{ z^{ t A } - 1 }{z} = (-1)^{ n } f_{ t } (z) \ , \]
so that
\begin{equation}\label{open} L ( {\cal P} , t ) = (-1)^{ n } \ \mbox{Res} \left( f_{ t } (z) , z = \infty \right) . \end{equation}
The finite poles of $f_{ t } $ are at 0 (with residue -1), 1, and the roots of unity in $\Omega$ as  before. This gives us

{\bf Theorem 1(b).} 
\[ L ( {\cal P} , t ) = (-1)^{ n } \left( 1 - \mbox{Res} \left( f_{ t } (z) , z=1 \right) - \sum_{ \lambda \in \Omega } \mbox{Res} \left( f_{ t } (z) , z=\lambda \right) \right) \ . \] 

As an immediate consequence we get the remarkable

{\bf Corollary 3 (Ehrhart-Macdonald Reciprocity Law)}. 
\[ L ( {\cal P} , -t ) = (-1)^{ n } L ( \overline{\cal P} , t ) \ . \]

This result was conjectured (again, for a general lattice polytope) by Ehrhart (\cite{ehrhart}), 
and later proved by Macdonald (\cite{macdonald}), McMullen (\cite{mcmullen}), and Stanley (\cite{stanley}).

\vspace{1cm}

{\bf 3. The Ehrhart coefficients}

With a small modification of $ f_{ t } (z) $, we can actually derive a formula for each coefficient of the Ehrhart
polynomial 
\[ L ( \overline{\cal P} , t ) = c_{n} t^{n} + \dots + c_{0} \ . \]
Consider the function
\begin{eqnarray*} g_{ k } (z) &:=& \frac{ \left( z^{ - t A } - 1 \right)^{k} }{ \left( 1 - z^{ A_{ 1 } }  \right) \left( 1 - z^{ A_{ 2 } }  \right) \cdots \left( 1 - z^{ A_{ n } }  \right) \left( 1 - z \right) z } \\
      &=& \frac{ \sum_{j=0}^{k} {k \choose j} z^{ - t A (k-j) } (-1)^{j} }{ \left( 1 - z^{ A_{ 1 } }  \right) \left( 1 - z^{ A_{ 2 } }  \right) \cdots \left( 1 - z^{ A_{ n } }  \right) \left( 1 - z \right) z } . \end{eqnarray*} 
If we insert $ - \sum_{j=0}^{k} {k \choose j} (-1)^{j} = 0 $ in the numerator, this becomes
\begin{eqnarray*} g_{ k } (z) &=& \sum_{j=0}^{k} {k \choose j} (-1)^{j} \frac{ z^{ - t (k-j) A } - 1}{ \left( 1 - z^{ A_{ 1 } }  \right) \left( 1 - z^{ A_{ 2 } }  \right) \cdots \left( 1 - z^{ A_{ n } }  \right) \left( 1 - z \right) z } \\
                              &=& \sum_{j=0}^{k-1} {k \choose j} (-1)^{j} f_{ -t(k-j) } (z) \ . \end{eqnarray*} 
Recall that (\ref{closure}) gave us $ L ( \overline{\cal P} , t ) = \mbox{Res} \left( f_{ -t } (z) , z=0 \right) + 1 $. 
Using this relation, we obtain
\begin{eqnarray*} &\mbox{}& \mbox{Res} \left( g_{ k } (z) , z=0 \right) = \sum_{j=0}^{k-1} {k \choose j} (-1)^{j} \ \mbox{Res} \left( f_{ -t(k-j) } (z) , z=0 \right) \\
      &\mbox{}& \quad = \sum_{j=0}^{k-1} {k \choose j} (-1)^{j} \Bigl( L \left( \overline{\cal P} , (k-j)t \right) - 1 \Bigr) \\
      &\mbox{}& \quad = \sum_{j=0}^{k-1} {k \choose j} (-1)^{j} L \left( \overline{\cal P} , (k-j)t \right) \ + \ (-1)^{k} \ . \end{eqnarray*} 
We claim that this polynomial has no terms with exponent smaller than $k$:

{\bf Lemma 4.} {\it Suppose $ L ( \overline{\cal P} , t ) = c_{n} t^{n} + \dots + c_{0} $. Then for $ 1 \leq k \leq n $ }
\begin{equation}\label{g_k} \mbox{Res} \left( g_{ k } (z) , z=0 \right) = k! \sum_{m=k}^{n} S(m,k) \ c_{m} \ t^{m} \ , \end{equation}
{\it where $S(m,k)$ denotes the Stirling number of the second kind.}

{\it Proof.} 
Suppose
\begin{equation}\label{b_km} \sum_{j=0}^{k-1} {k \choose j} (-1)^{j} L \left( \overline{\cal P} , (k-j)t \right) = \sum_{m=0}^{n} b_{k,m} t^{m} \ , \end{equation}
so that for $m>0$
\[ b_{k,m} = \sum_{j=0}^{k-1} {k \choose j} (-1)^{j} c_{m} (k-j)^{m} = c_{m} \sum_{j=0}^{k} {k \choose j} (-1)^{k-j} j^{m} \ . \]
The Stirling number of the second kind $S(m,k)$ is the number of partitions of an $m$-set into 
$k$ blocks. The reason we are interested in these numbers is the identity (\cite{stanley})
\[ S(m,k) = \frac{ 1 }{ k! } \sum_{j=0}^{k} {k \choose j} (-1)^{k-j} j^{m} , \]
so that $ b_{k,m} = c_{m} \ k! \ S(m,k) $ for $m>0$. Some of the elementary properties of $ S(m,k) $ are (\cite{stanley})
\begin{eqnarray} &\mbox{}& S(m,k) = 0 \quad \mbox{ if } k>m \label{zero} \\
                 &\mbox{}& S(m,1) = 1 \label{S_m1} \\
                 &\mbox{}& S(m,m) = 1 \label{S_mm} \\
                 &\mbox{}& S(m,k) = k \ S(n-1,k) + S(n-1, k-1) \ . \nonumber  \end{eqnarray} 
By (\ref{zero}), we conclude that $b_{k,m} = 0$ for $1 \leq m < k$. 
The constant term in (\ref{b_km}) is
\[ b_{k,0} = \sum_{j=0}^{k-1} {k \choose j} (-1)^{j} c_{0} = - c_{0} (-1)^{ k } . \]
Since $c_{0} = 1$ for our tetrahedron (in fact, $c_{0} = 1$ for any convex polytope \cite{ehrhart}), 
(\ref{g_k}) follows. \hfill {} $\Box$

The other poles of $g_{ k }$ are at 1 and the roots of unity in 
\[ \Omega_{k} := \left\{ z \in \C \setminus \{ 1 \} : z^{ \frac{A}{ a_{j_{1}} \cdots a_{j_{k+1}}} } = 1 , 1 \leq j_{1} < j_{2} < \dots < j_{k+1} \leq n \right\} . \]
Note that as $k$ gets larger, $\Omega_{k}$ gets smaller. That is, we have fewer residues to 
consider. This is consistent with the notion that the computational complexity increases with 
each additional coefficient, that is, the computation of $c_{k}$ is more complicated than that of $c_{k+1}$. 
Using the residue theorem, we can rewrite (\ref{g_k}) as

{\bf Theorem 5.} {\it Suppose $ L ( \overline{\cal P} , t ) = c_{n} t^{n} + \dots + c_{0} $. Then for $ 1 \leq k \leq n $ }
\[ \sum_{m=k}^{n} S(m,k) \ c_{m} \ t^{m} = \frac{-1}{k!} \left( \mbox{Res} \left( g_{ k } (z) , z=1 \right) + \sum_{\lambda \in \Omega_{k}} \mbox{Res} \left( g_{ k } (z) , z=\lambda \right) \right) . \]

{\it Remarks.} 1. For $k=1$, we get with (\ref{S_m1}) a reformulation of Theorem 1(a). 

2. The coefficients of $ L ( {\cal P} , t ) $ are the same as those of $ L ( \overline{\cal P} , t ) $, up to the sign: By Corollary 3, 
$ L ( {\cal P} , t ) = c_{n} t^{n} - c_{n-1} t^{n-1} + \dots + (-1)^{n} c_{0}  $. 

3. Res($g_{ k } (z) , z=1$) can be computed as easily as before, the slightly more difficult task is to
get the residues at the roots of unity (see also remark 2 following Theorem 1(a)). However, with increasing $k$, we have to consider fewer of
them, so that there is less to calculate. If we want to compute the Ehrhart coefficient $c_{m}$,
we only have to consider the roots of unity in $\Omega_{m}$. We can make this more precise:
With (\ref{S_mm}), we obtain

{\bf Corollary 6.} {\it For $m>0$, $c_{m}$ is the coefficient of $ t^{ m } $ of  }
\[ \frac{ -1 }{ m! } \left( \mbox{Res} \left( g_{ m } (z) , z=1 \right) + \sum_{\lambda \in \Omega_{m}} \mbox{Res} \left( g_{ m } (z) , z=\lambda \right)  \right)  . \]

\vspace{1cm}

{\bf 4. An Example}

As an application, we will compute the first nontrivial Ehrhart coefficient $c_{n-2}$ for the 
$n$-dimensional tetrahedron ${\cal P}$ ($n \geq 3$) under the additional assumption that  
$a_{ 1 } , \dots , a_{ n } $ are pairwise relatively prime integers $\geq 2$. 
This case was first explored by Pommersheim (\cite{pomm}). 

{\bf Theorem 7.} {\it Under the above assumptions, }
\[ c_{n-2} = \frac{1}{(n-2)!} \ \Bigl( C_{n} - \s (A_{1},a_{1}) - \dots - \s (A_{n},a_{n}) \Bigr) ,  \]
{\it where $\s (a,b)$ denotes the Dedekind sum, and}
\[ C_{n} := \frac{1}{4} \left( n + A_{1,2} + \dots + A_{n-1,n} \right) + \frac{1}{12} \left( \frac{1}{A} + \frac{ A_{1} }{ a_{1} } + \dots + \frac{ A_{n} }{ a_{n} }  \right) \ . \]
Here $ A_{j,k} $ denotes $ a_{1} \cdots {\hat a}_{j} \cdots {\hat a}_{k} \cdots a_{n} $. 

{\it Proof.} We have to consider
\[ g_{ n-2 } (z) = \frac{ \left( z^{ - t A } - 1 \right)^{n-2} }{ \left( 1 - z^{ A_{ 1 } }  \right) \left( 1 - z^{ A_{ 2 } }  \right) \cdots \left( 1 - z^{ A_{ n } }  \right) \left( 1 - z \right) z } \]
Because $a_{ 1 } , \dots , a_{ n } $ are pairwise relatively prime, 
$g_{ n-2 }$ has {\it simple} poles at all the $ a_{1} , \dots , a_{n} $th roots of unity. Let 
$\lambda^{a_{1}} = 1 \not= \lambda$. Then 
\begin{eqnarray*} &\mbox{}& \mbox{Res} \left( g_{ n-2 } (z) , z = \lambda \right) = \\
                  &\mbox{}& \quad = \frac{1}{ \left( 1 - \lambda^{A_{1}} \right) \left( 1 - \lambda \right) \lambda } \ \mbox{Res} \left( \frac{ \left( z^{ - t A } - 1 \right)^{n-2} }{\left( 1 - z^{ A_{ 2 } }  \right) \cdots \left( 1 - z^{ A_{ n } }  \right)} , z = \lambda \right) . \end{eqnarray*}
Similar to the methods used in the first chapter to arrive at Corollary 2, we make a 
change of variables $ z = w^{1/a_{1}} = \exp \left( \frac{1}{a_{1}} \log w \right) $, where we choose 
a suitable branch of the logarithm such that $ \exp \left( \frac{1}{a_{1}} \log (1) \right) = \lambda$. 
We thus obtain
\begin{eqnarray*} &\mbox{}& \mbox{Res} \left( g_{ n-2 } (z) , z = \lambda \right) = \\
                  &\mbox{}& \quad = \frac{1}{ \left( 1 - \lambda^{A_{1}} \right) \left( 1 - \lambda \right) \lambda } \frac{\lambda}{a_{1}} \ \mbox{Res} \left( \frac{ \left( w^{ - t B } - 1 \right)^{n-2} }{\left( 1 - w^{ B_{ 2 } }  \right) \cdots \left( 1 - w^{ B_{ n } }  \right)} , w = 1 \right) , \end{eqnarray*}
where $B := a_{2} \cdots a_{n} , \ B_{k} := a_{2} \cdots {\hat a}_{k} \cdots a_{n} $. We claim that
\[ \mbox{Res} \left( \frac{ \left( z^{ - t B } - 1 \right)^{n-2} }{\left( 1 - z^{ B_{ 2 } }  \right) \cdots \left( 1 - z^{ B_{ n } }  \right)} , z = 1 \right) = - t^{n-2} .  \]
To prove this, first note that 
\[ \left( z^{ - t B } - 1 \right)^{n-2} = (-tB)^{n-2} (z-1)^{ n-2 }  + O \left( (z-1)^{ n-1 } \right) .  \]
Now for $m \in \N$, 
\[ \mbox{Res} \left( \frac{1}{1 - z^{m}} , z=1 \right) = \lim_{z \to 1} \frac{z-1}{1 - z^{m}} = - \frac{1}{m} \ . \]
Putting all of this together, we obtain
\begin{eqnarray*} &\mbox{}& \mbox{Res} \left( \frac{ \left( z^{ - t B } - 1 \right)^{n-2} }{\left( 1 - z^{ B_{ 2 } }  \right) \cdots \left( 1 - z^{ B_{ n } }  \right)} , z = 1 \right) = \frac{ (-tB)^{n-2} }{ (-B_{2}) \cdots (-B_{n}) } \\
                  &\mbox{}& \quad = - \frac{ t^{n-2} a_{2}^{n-2} \cdots a_{n}^{n-2} }{a_{2}^{n-2} \cdots a_{n}^{n-2}} = - t^{n-2} , \end{eqnarray*} 
as desired. Therefore 
\[ \mbox{Res} \left( g_{ n-2 } (z) , z = \lambda \right) = \frac{ - t^{n-2} }{ a_{1} \left( 1 - \lambda^{A_{1}} \right) \left( 1 - \lambda \right) } \ . \]
Adding up all the $a_{1}$'th roots of unity $\not= 1$, we get 
\begin{eqnarray*} &\mbox{}& \sum_{ \lambda^{ a_{1} } = 1 \not= \lambda } \mbox{Res} \left(g_{ n-2 } (z) , z = \lambda \right) = \frac{- t^{n-2}  }{ a_{1} } \sum_{ \lambda^{ a_{1} } = 1 \not= \lambda } \frac{ 1 }{ \left( 1 - \lambda^{A_{1}} \right) \left( 1 - \lambda \right) } \\
                  &\mbox{}& \quad = \frac{- t^{n-2}  }{ a_{1} } \sum_{ k=1 }^{a_{1} - 1 } \frac{ 1 }{ \left( 1 - \xi^{ k A_{1} }  \right) \left( 1 - \xi^{ k }  \right)  } , \end{eqnarray*} 
where $\xi$ is a primitive $a_{1}$'th root of unity. This finite sum is practically a Dedekind sum: 
\begin{eqnarray*} &\mbox{}& \frac{ 1 }{ a_{1} } \sum_{ k=1 }^{ a_{1}-1 } \frac{ 1 }{ \left( 1 - \xi^{ kA_{1} }  \right) \left( 1 - \xi^{ k }  \right)  } = \frac{ 1 }{ 4a_{1} } \sum_{ k=1 }^{ a_{1}-1 } \left( 1 + \frac{ 1 + \xi^{ kA_{1} }  }{ 1 - \xi^{ kA_{1} }  }  \right) \left( 1 + \frac{ 1 + \xi^{ k }  }{ 1 - \xi^{ k }  }  \right)   \\
                  &\mbox{}& \quad = \frac{ 1 }{ 4a_{1} } (a_{1}-1) - \frac{ i }{ 4a_{1} } \sum_{ k=1 }^{ a_{1}-1 } \left( \cot \frac{ \pi k A_{1} }{ a_{1} } + \cot \frac{ \pi k }{ a_{1} } \right) \\
                  &\mbox{}& \qquad - \frac{ 1 }{ 4a_{1} } \sum_{ k=1 }^{ a_{1}-1 } \cot \frac{ \pi k A_{1} }{ a_{1} } \cot \frac{ \pi k }{ a_{1} } \\
                  &\mbox{}& \quad = \frac{ 1 }{ 4 } - \frac{ 1 }{ 4a_{1} } - \s (A_{1},a_{1}) . \end{eqnarray*} 
The imaginary terms disappear here, since the sum on the left hand side and $\s(A_{1},a_{1})$ 
are rational: Both are elements of the cyclotomic field of $a_{1}$'th roots of unity, and 
invariant under all Galois transformations of this field. 

Hence we obtain
\[ \sum_{ \lambda^{ a_{1} } = 1 \not= \lambda } \mbox{Res} \left(g_{ n-2 } (z) , z = \lambda \right) = - t^{n-2} \left( \frac{ 1 }{ 4 } - \frac{ 1 }{ 4a_{1} } - \s (A_{1},a_{1}) \right) . \]
We get similar expressions for the residues at the other roots of unity, so that Corollary 6 gives us for $n \geq 3$
\begin{eqnarray}\label{c_n-2} &\mbox{}& c_{n-2} = \frac{1}{(n-2)!} \biggl( \frac{ n }{ 4 } - C - \frac{ 1 }{ 4 } \left( \frac{ 1 }{ a_{1}  } + \dots + \frac{ 1 }{ a_{n}  } \right) \biggr. \nonumber \\
                  &\mbox{}& \qquad \qquad \qquad \biggl. - \s (A_{1},a_{1}) - \dots - \s (A_{n},a_{n}) \biggr) ,  \end{eqnarray}
where $C$ is the coefficient of $ t^{ n-2 } $ of $ \mbox{\rm Res} \left( g_{n-2} (z) , z = 1 \right) $.
We can actually obtain a closed form for $C$: As before, 
\begin{eqnarray*} &\mbox{}& \mbox{Res} \left( g_{ n-2 } (z) , z=1 \right) = \mbox{Res} \left( e^{ z }  g_{ n-2 } (e^{ z } ) , z=0 \right)  \\
                  &\mbox{}& \quad = \mbox{Res} \left( \frac{ \left( e^{ -t A z } - 1 \right)^{n-2} }{ \left( 1 - e^{ A_{ 1 } z } \right) \left( 1 - e^{ A_{ 2 } z } \right) \cdots \left( 1 - e^{ A_{ n } z } \right) \left( 1 - e^{ z } \right) } , z=0 \right) . \end{eqnarray*} 
Now with 
\[ \left( e^{ -t A z } - 1 \right)^{n-2} = \left( -t A z \right)^{n-2} + O \left( \left( t z \right)^{n-1}  \right) \]
and 
\[ \frac{ 1 }{ 1 - e^{ z }  } = - z^{ -1 } + \frac{ 1 }{ 2 } - \frac{ 1 }{ 12 } z + O\left( z^{ 3 } \right) , \]
the coefficient of $ t^{ n-2 } $ of $ \mbox{\rm Res} \left( g_{n-2} (z) , z = 1 \right) $ turns out to be
\begin{eqnarray*} &\mbox{}& C = \left( -A \right)^{n-2} \left[ \frac{1}{12} \left( \frac{ (-1)^{ n+1 } }{ A_{1} \cdots A_{n} } + \frac{ (-1)^{ n+1 } A_{1}  }{ A_{2} \cdots A_{n} } + \dots + \frac{ (-1)^{ n+1 } A_{n}  }{ A_{1} \cdots A_{n-1} }  \right) + \right. \\
                  &\mbox{}& \qquad \qquad \qquad \qquad + \frac{1}{4} \left( \frac{ (-1)^{ n-1 } }{ A_{2} \cdots A_{n} } + \dots + \frac{ (-1)^{ n-1 } }{ A_{1} \cdots A_{n-1} } \ + \right. \\
                  &\mbox{}& \qquad \qquad \qquad \qquad \left. \left.  + \frac{ (-1)^{ n-1 } }{ A_{3} \cdots A_{n} } + \frac{ (-1)^{ n-1 } }{ A_{2} A_{4} \cdots A_{n} } + \dots + \frac{ (-1)^{ n-1 } }{ A_{1} \cdots A_{n-2} }  \right) \right] \\ 
                  &\mbox{}& \quad = - \frac{1}{12} \left( \frac{1}{A} + \frac{ A_{1} }{ a_{1} } + \dots + \frac{ A_{n} }{ a_{n} }  \right) \\
                  &\mbox{}& \qquad \qquad - \frac{1}{4} \left( \frac{1}{ a_{1} } + \dots + \frac{1}{ a_{n} } + A_{1,2} + \dots + A_{n-1,n} \right) . \end{eqnarray*} 
Substituting this into (\ref{c_n-2}) yields the statement. \hfill {} $\Box$

The other Ehrhart-coefficients for this tetrahedron can be derived in a similar fashion, 
although the computation gets more and more complicated, as noted in section 3. 

\vspace{1cm}

{\bf 5. General lattice polytopes}

Any convex lattice polytope (that is, a convex polytope whose vertices are on the lattice $\Z^{n}$) 
can be described by a finite number of inequalities over the integers. In other words, a 
convex lattice polytope ${\cal P} $ is an intersection of finitely many half-spaces. Translation does 
not change the lattice point count, so we can assume that the points in the polytope have 
positive coordinates and apply the ideas of the previous sections to ${\cal P} $. 
Suppose the closure of the dilated polytope $ \overline{t {\cal P}} $ is given by the $n+q$ inequalities
\begin{eqnarray}\label{inequ} &\mbox{}& x_{ 1 } , \dots , x_{ n } \geq 0 \nonumber \\ \nonumber \\ 
                  &\mbox{}& \frac{ x_{ 1 }  }{ a_{ 11 }  } + \dots + \frac{ x_{ n }  }{ a_{ 1n }  } \leq t \\
                  &\mbox{}& \qquad \vdots \nonumber \\
                  &\mbox{}& \frac{ x_{ 1 }  }{ a_{ q1 }  } + \dots + \frac{ x_{ n }  }{ a_{ qn }  } \leq t \ , \nonumber \end{eqnarray} 
with $a_{ jk } \in \Z$. Then by introducing a similar notation as before, 
\[ P_{k} := a_{k1} \cdots a_{kn} , \quad A_{jk} := a_{j1} \cdots {\hat a}_{jk} \cdots a_{jn} , \]
we can define a matrix
\[ M = \Bigl( A_{jk} \Bigr)_{ j=1..q \atop k=1..n } . \]
Let $C_{j}$ the $j$'th column, and $R_{k}$ the $k$'th row of $M$. Then we can rewrite the 
last $q$ inequalities of (\ref{inequ}) as 
\begin{eqnarray}\label{inequ2} &\mbox{}& R_{ 1 } \cdot x \leq t P_{1} \nonumber \\
                  &\mbox{}& \qquad \vdots \\
                  &\mbox{}& R_{ q } \cdot x \leq t P_{q} , \nonumber \end{eqnarray} 
where $ x = ( x_{ 1 } , \dots , x_{ n } ) $ and $\cdot$ denotes the usual scalar product. 
Now consider the function
\[ f(z) = f \left( z_{ 1 } , \dots , z_{ q }  \right) = \frac{ z_{ 1 }^{ -t P_{1} - 1 } \cdots z_{ q }^{ -t P_{q} - 1 }    }{ \left( 1 - z^{ C_{ 1 }  } \right) \cdots \left( 1 - z^{ C_{ q }  } \right) \left( 1 - z_{ 1 }  \right) \cdots \left( 1 - z_{ q }  \right)  } . \]
Here we used the standard multinomial notation $ z^{v} := z_{1}^{ v_{1} } \cdots z_{q}^{ v_{q} } $. 
We will integrate $f$ with respect to each variable over a circle with small radius: 
\begin{equation}\label{integral} \int_{ \left| z_{ 1 }  \right| = \epsilon_{ 1 } } \cdots \int_{ \left| z_{ q }  \right| = \epsilon_{ q } } f \left( z_{ 1 } , \dots , z_{ q }  \right) \ dz_{ q } \cdots dz_{ 1 }  \end{equation}
Here, $ 0 < \epsilon_{0}, \dots , \epsilon_{q} < 1 $ are chosen such that we can expand all the $ \frac{ 1 }{ 1 - z^{ C_{k} } } $ into power series 
about 0. To ensure the existence of $ \epsilon_{0}, \dots , \epsilon_{q} $, we may, if necessary, add an additional inequality 
$ x_{1} + \dots + x_{n} \leq t P_{0} $ for a suitable large $ P_{0} $. This is always possible, since $ {\cal P} $ is bounded. 

Since the integral over one variable will give the respective residue at 0, we 
can integrate with respect to one variable at a time. When $f$ is expanded into its Laurent 
series about 0, each term has the form
\[ z_{ 1 }^{ m \cdot R_{ 1 } + r_{ 1 } - t P_{ 1 } - 1 } \cdots z_{ q }^{ m \cdot R_{ q } + r_{ q } - t P_{ q } - 1 } ,  \]
where $ m := ( m_{ 1 } , \dots , m_{ n } ) $, and $ m_{ 1 } , \dots , m_{ n } , r_{ 1 } , \dots , r_{ q } $ are nonnegative integers. Thus, in 
the integral (\ref{integral}), this term will give a contribution 
precisely if $ m $ satisfies the inequalities (\ref{inequ2}). In other words, we have proved

{\bf Theorem 8.} 
\[ L ( \overline{\cal P} , t ) = \frac{ 1 }{ (2 \pi i)^{ q }  } \int_{ \left| z_{ 1 }  \right| = \epsilon_{ 1 } } \cdots \int_{ \left| z_{ q }  \right| = \epsilon_{ q } } f \left( z_{ 1 } , \dots , z_{ q }  \right) \ dz_{ q } \cdots dz_{ 1 } \ . \]

On the other hand, we can compute this integral by computing the residues $ \not= 0$ as before, 
for each variable at a time. The (finite) sum over these residues will be equal, by Theorem 10, 
to the number of lattice points in the closure of our polytope $\overline{t {\cal P}}$. We will explore this 
method and its results in a future paper.

\vspace{1cm}

{\bf 6. Further remarks}

1. The ideas in this paper can also be used in other ways. First, note that by omitting the 
$z^{ \pm t A } - 1$ factor, we essentially get trivial residues at 
$\infty$ and 0. In view of Corollary 6, this corresponds to obtaining an expression for the constant Ehrhart 
coefficient and setting it equal to 1 (\cite{beck}). This means we are not getting a lattice point count; however, 
the other factors will still give us Dedekind-like sums, as they appear in the formulae. By using 
the residue theorem as before, we thus now get relations between these terms, similar in spirit 
to the various reciprocity laws for classical (\cite{grosswald}) and higher dimensional 
Dedekind sums (\cite{zagier}). 

2. Finally, we can extend these methods to $n$-dimensional {\it rational} polytopes 
(whose vertices have now rational coordinates). 
For example, the number of lattice points in the dilated rational tetrahedron
\[ \left\{ (x_{ 1 } , \dots , x_{ n }) \in \R^{ n }  : \sum_{ k=1 }^{ n } x_{ k } a_{ k } < t \ \mbox{and all} \ x_{ k } > 0  \right\} \]
can be computed with similar methods. The details will appear in a future 
paper (\cite{rational}). 

\vspace{1cm}

{\bf Acknowledgements}

The author would like to thank Boris Datskovsky, Ira Gessel, Marvin Knopp, Gerardo Mendoza, and especially Sinai Robins for helpful suggestions and encouragement.

\vspace{1cm}

\footnotesize
\nocite{*}
\addcontentsline{toc}{subsubsection}{References}
\bibliography{thesis}
\bibliographystyle{alpha}

\vspace{1cm}

 \sc Department of Mathematical Sciences\\
 State University of New York\\
 Binghamton, NY 13902-6000\\
 {\tt matthias@math.binghamton.edu} 

\end{document}